\documentclass[12pt]{amsart}
\usepackage{amsmath, amsthm, amssymb}
\usepackage{fullpage}
\usepackage{color}
\usepackage{hyperref}
\usepackage{soul}
\usepackage{enumitem}
\usepackage{cancel}
\usepackage[normalem]{ulem}

\newtheorem{theorem}{Theorem}
\newtheorem{lemma}{Lemma}
\newtheorem{proposition}[lemma]{Proposition}

\newtheorem{definition}[lemma]{Definition}
\newtheorem{remark}[lemma]{Remark}
\newtheorem{conjecture}{Conjecture}
\numberwithin{lemma}{section}

\DeclareMathOperator{\sech}{sech}

\numberwithin{equation}{section}

\newcommand{\R}{{\mathbb R}}

\newcommand{\Z}{{\mathbb Z}}

\renewcommand{\R}{\mathbb R}

\newcommand{\bM}{\mathbf M}
\newcommand{\bP}{\mathbf P}

\newcommand{\bI}{\mathbf I}
\newcommand{\bJ}{\mathbf J}
\newcommand{\bK}{\mathbf K}

\newcommand{\bu}{{\bar u}}

\newcommand{\la}{\langle}
\newcommand{\ra}{\rangle}

\newcommand{\ms}{M^\sharp}
\newcommand{\ps}{P^\sharp}
\newcommand{\calR}{\mathcal{R}}

\newcommand{\xih}{{\delta \xi^{\text{hi}}}}
\newcommand{\xim}{{\delta \xi^{\text{med}}}}

\newcommand{\tDelta}{{\tilde \Delta}}

\begin{document}

\title{Long time solutions for 1D cubic dispersive equations, Part II: the focusing case}

\author{Mihaela Ifrim}
\address{Department of Mathematics, University of Wisconsin, Madison}
\email{ifrim@wisc.edu}

\author{ Daniel Tataru}
\address{Department of Mathematics, University of California at Berkeley}
\email{tataru@math.berkeley.edu}

\begin{abstract}
This article is concerned with one dimensional dispersive flows with cubic nonlinearities on the real line.  In a very recent work, the authors have introduced a broad conjecture for such flows, asserting that
in the defocusing case, small initial data yields global, scattering solutions. Then this conjecture was proved in the case of a Schr\"odinger dispersion relation. In terms of scattering, our global solutions were proved to satisfy both global $L^6$ Strichartz estimates and bilinear $L^2$ bounds. Notably, no localization assumption is made on the initial data.

In this article we consider the focusing scenario. There potentially one may have small solitons, so one cannot hope to have global scattering solutions in general. Instead, we look for long time solutions, and ask what is the time-scale on which the solutions exist and satisfy good dispersive estimates.
Our main result, which also applies in the case of the Schr\"odinger dispersion relation, asserts that for initial data of size $\epsilon$, the solutions exist on the time-scale $\epsilon^{-8}$, and satisfy the 
desired  $L^6$ Strichartz estimates and bilinear $L^2$ bounds on the time-scale $\epsilon^{-6}$. 
To the best of our knowledge, this is the first result to reach such a threshold.

\end{abstract}

\maketitle

\setcounter{tocdepth}{1}
\tableofcontents


\section{Introduction}
The question of obtaining long time solutions for one dimensional dispersive flows with quadratic/cubic nonlinearities has attracted a lot of attention in recent years. One can distinguish two different but closely related types of results
that have emerged, as well as several successful approaches.

On one hand, \emph{normal form methods} have been developed 
in order to extend the lifespan of solutions, beginning with 
\cite{Shatah} in the late '80's. Somewhat later, around 2000, the \emph{ I-method}, introduced in \cite{I-method} brought forth the idea of constructing  better almost conserved quantities. These two ideas serve well in the study of semilinear flows, where it was later understood that they are connected ~\cite{Bourgain-nf}.

Neither of these techniques can be directly applied to quasilinear problems. Addressing this problem, it was discovered in the work of the authors and collaborators \cite{BH},
\cite{IT-g} that one can adapt the normal form method 
to quasilinear problems by constructing energies which 
simultaneously capture both the quasilinear and the normal form structures. This idea was called the \emph{modified energy method}, and can also be seen in some way as a quasilinear adaptation of the I-method. 
An alternate approach, also in the quasilinear setting, is provided by the \emph{ flow method } of 
\cite{hi}, where a better normal form transformation
is constructed using a well chosen auxiliary flow.

On the other hand, the further goal of obtaining
scattering, global in time solutions for one dimensional dispersive flows with quadratic/cubic nonlinearities  has also been extensively studied in the last two decades  for a number of models, under the assumption that the initial data is both \emph{small} and \emph{localized};  without being exhaustive, see for instance \cite{HN,HN1,LS,KP,IT-NLS}.
The nonlinearities in these models are primarily cubic, though the analysis has also been extended via normal form and modified energy methods to problems which also have nonresonant quadratic interactions; several such examples are  \cite{AD,IT-g,D,IT-c,LLS}, see also further references therein,
as well as the authors' expository paper \cite{IT-wp}.

If instead one considers initial data which is just \emph{small}, without any localization assumption, then the problem becomes much more difficult, because this allows for far stronger nonlinear interactions over long time-scales. One also needs to distinguish between the focusing 
and the defocusing problems. In a recent paper \cite{IT-noloc}, the authors have introduced a broad global well-posedness (GWP) conjecture, which applies to both semilinear and quasilinear problems:

\begin{conjecture}[Non-localized data defocusing GWP conjecture]\label{c:nld}
One dimensional dispersive flows on the real line with cubic defocusing nonlinearities and small initial data have global in time, scattering solutions.
\end{conjecture}

The main result of \cite{IT-noloc} asserts that this conjecture is true under suitable assumptions, most notably that the dispersion relation is the Schr\"odinger dispersion relation. 
That was the first global in time well-posedness result of this type. Notably, scattering 
here is interpreted in a weak sense, to mean that the solution satisfies global $L^6$ Strichartz estimates and bilinear $L^2$ bounds. This is because of the strong nonlinear effects, which preclude
any kind of classical scattering. The precise result is stated later in Theorem~\ref{t:defocusing}.

\medskip

Our interest in this article is instead in the focusing case of the same problem. Since 1D focusing dispersive problems typically admit small solitons, a global result as stated in the above conjecture
simply cannot hold. Even if global solutions exist (as it is the case for instance for the cubic focusing NLS problem) the presence of solitons will defeat any kind of global decay estimates.
For this reason,  in the focusing case we will rethink the problem as a question about the lifespan of solutions with small initial data. Precisely, if the initial data has size $\epsilon$ when measured 
in a suitable $H^s$ Sobolev norm, what can be said about the lifespan of the solutions as a function of $\epsilon$ ?

Following the lead of our earlier paper, we begin by formulating the focusing counterpart of the previous conjecture. Then we will prove that the conjecture is true under suitable assumptions. Our main conjecture is as follows:

\begin{conjecture}[Non-localized data focusing conjecture]
One dimensional dispersive flows on the real line with cubic nonlinearities and small initial data of size $\epsilon$ have solutions which remain of comparable size at least on an $\epsilon^{-8}$ time-scale.
\end{conjecture}

The main result of this paper, see Theorem~\ref{t:focusing} below, asserts that this conjecture is valid under the additional assumption that the dispersion relation is of Schr\"odinger type. As part of this result, we also prove that our long time solutions  satisfy both $L^6$ Strichartz estimates and bilinear $L^2$ bounds on suitable time-scales. This is akin to our earlier work on the defocusing case, but with the  difference that in the defocusing case such estimates were proved globally in time.

For reference purposes, we note some intermediate lifespan thresholds
which can be reached with methods which were developed earlier:

\begin{itemize}
    \item A cubic lifespan $\epsilon^{-2}$ can be reached using direct energy estimates, using only the fact that the nonlinearity is cubic.
    
    \item A quintic lifespan $\epsilon^{-4}$ can be reached by more accurate energy estimates using a quartic energy correction. This requires  
    the cubic nonlinearity to be conservative, as defined later in the introduction.
    
    \item In the case of a perturbative nonlinearity, a quintic  lifespan can also be 
    obtained by directly using Strichartz estimates.
\end{itemize}

Heuristically, if it were possible to directly combine the last two ideas above, that would lead to the $\epsilon^{-8}$ threshold in the present paper. However, the price to achieve that would be very steep, as one would need to assume both a perturbative nonlinearity and high regularity. By comparison, our new result provides a  much more robust approach, which is both of a nonperturbative nature and far more efficient in terms of the regularity requirements.
We show here how dispersive and normal form tools can be combined  very efficiently 
in order to make substantial gains. 

 In the present paper we  aim for a reasonably simple setting, where our model problem is borderline semilinear, and where we prove results we expect to be optimal. This should also serve as a  baseline for further developments. In particular, we believe that our methods can be also applied in non-perturbative, quasilinear settings.

\subsection{Cubic NLS problems in one space dimension}
The fundamental model for one-dimensional dispersive flows 
with cubic nonlinearity  in one space 
dimension is the cubic nonlinear Schr\"odinger (NLS) flow, 
\begin{equation}\label{nls3}
\left\{
\begin{aligned}
&i u_t + u_{xx} = \pm 2u |u|^2 \\
&u(0) = u_0,
\end{aligned}
\right.
\end{equation}
with $u:\mathbb{R}\times \mathbb{R}\rightarrow \mathbb{C}$. This comes in a defocusing (+) and a focusing (-) version.

The above cubic NLS flows are globally well-posed in $L^2$ both in the focusing and in the defocusing case, though the global behavior differs in the two cases. Both of these model problems are completely integrable,
and one may study their global behavior using inverse scattering tools \cite{IST}, \cite{IST-focusing}.

In the defocusing case, the inverse scattering approach  allows one to treat
the case of localized data, and show that global solutions scatter at infinity, see for instance \cite{IST}. This can also be proved in a more robust way, without using inverse scattering, under the assumption that the initial data is small and localized, see \cite{IT-NLS} and references therein.
Much less is known in terms of scattering for nonlocalized $L^2$ data. However, if more regularity is assumed for the data, then we have the following estimate due to Planchon-Vega~\cite{PV}, see also the work of Colliander-Grillakis-Tzirakis~\cite{MR2527809}:
\begin{equation}\label{PV-est}
\|u\|_{L^6_{t,x}}^6 + \| \partial_x |u|^2\|_{L^2_{t,x}}^2 \lesssim \| u_0\|_{L^2_x}^3 \|u_0\|_{H^1_x}.
\end{equation}
This allows one to estimate the $L^6$ Strichartz norm of the solution, i.e. to prove some type of scattering or dispersive decay. This estimate was improved and extended to $L^2$ solutions as a corollary of the results in our previous paper \cite{IT-noloc}. Precisely, we have 
\begin{equation}\label{PV-est-new}
\|u\|_{L^6_{t,x}}^2 + \| \partial_x |u|^2\|_{L^2_t(\dot H^{-\frac12}_x + c L^2_x)}  \lesssim \| u_0\|_{L^2_x}^2, \qquad c = \|u_0\|_{L^2_x}. 
\end{equation}

On the other hand, the focusing problem admits small solitons, so the solutions cannot in general scatter at infinity. If in addition the initial data is localized, then one expects the solution to resolve into a superposition of (finitely many) solitons, and a dispersive part; this is called \emph{the soliton resolution conjecture}, and is known to hold in a restrictive setting, via the method of inverse scattering, see e.g. ~\cite{IST-focusing}.

\bigskip

Our interest here is in focusing problems, but without any integrability
assumptions, and even without assuming any conservation laws. The model 
we consider is similar to the one in \cite{IT-noloc}, namely
\begin{equation}\label{nls}
\left\{
\begin{aligned}
&i u_t + u_{xx} = C(u,\bar u, u)\\
&u(0) = u_0,
\end{aligned}
\right.
\end{equation}
where $u$ is a complex valued function, $u:\mathbb{R}\times
\mathbb{R}\rightarrow \mathbb{C}$. Here $C$ is a trilinear translation
invariant form, whose symbol $c(\xi_1,\xi_2,\xi_3)$ can always be
assumed to be symmetric in $\xi_1,\xi_3$; see \cite{IT-noloc}
for an expanded discussion of  multilinear forms.
The arguments $u,\bar u$ and $u$ of $C$ are chosen so that
our equation \eqref{nls} has the phase rotation symmetry,
$u \to u e^{i\theta}$, as it is the case in many examples of interest.
The symbol $c(\xi_1,\xi_2,\xi_3)$ will be required to satisfy 
the following set of assumptions, which are similar to \cite{IT-noloc}:

\begin{enumerate}[label=(H\arabic*)]

\item Bounded and regular: 
\begin{equation}\label{c-smooth}
|\partial_\xi^\alpha c(\xi_1,\xi_2,\xi_3)| \leq c_\alpha,
\qquad \xi_1,\xi_2,\xi_3 \in \R,\,   \mbox{  for every  multi-index $\alpha$}.
\end{equation}

\item Conservative: 
\begin{equation}\label{c-conserv}
\Im c(\xi,\xi,\eta) = 0, \qquad \xi,\eta \in \R, \mbox { where } \Im z=\mbox{imaginary part of } z\in \mathbb{C}.
\end{equation}

\end{enumerate}

In addition to these two conditions, in \cite{IT-noloc} we have also employed
a defocusing assumption, namely

\begin{enumerate}[label=(H\arabic*),resume]
\item Defocusing: 
\begin{equation}\label{c-defocus}
c(\xi,\xi,\xi) \geq c_0 > 0, \qquad \xi \in \R  \mbox{ and } c_0 \in \mathbb{R^+}.
\end{equation}
\end{enumerate}

Here one might think that we should require the opposite, namely
\begin{enumerate}[label=(H\arabic*),resume]
\item Focusing: 
\begin{equation}\label{c-focus}
- c(\xi,\xi,\xi) \geq c_0 > 0, \qquad \xi \in \R  \mbox{ and } c_0 \in \mathbb{R^+}.
\end{equation}
\end{enumerate}
But as it turns out, no such assumption is needed here, as the result 
of this paper applies equally regardless of any sign condition; so we will simply drop it. 

Using the  same assumptions (H1), (H2) as in \cite{IT-noloc} is convenient here because it will allow us to reuse a good part of the analysis there,
up to the point where the defocusing assumption is needed.

Repeating a similar comment in \cite{IT-noloc}, one should view both our 
choices of the Schr\"odinger dispersion relation and the uniform bounds in (H1) not as fundamental, but rather as a balance between the generality of the result on one hand, and a streamlined exposition on the other hand.  This choice places 
our model in the semilinear class, but just barely so.

The simplest example of such a trilinear form $C$ is of course $C = \pm 1$, which corresponds to the classical one-dimensional cubic NLS problem.
But this problem has too much  structure, in particular it is  completely integrable, and also globally well-posed in $L^2$.

At the other end, both our use of the linear Schr\"odinger operator
and the boundedness condition (H1) are non-optimal, and we hope to relax
both of these restrictions in subsequent work.

\subsection{ The main result}

Our main result asserts that long time well-posedness holds for 
our problem for small $L^2$ data.  In addition, our solutions not only satisfy uniform $L^2$, but also space-time $L^6$ Strichartz estimates,
as well as bilinear $L^2$ bounds, on appropriate time-scales:

\begin{theorem}\label{t:focusing}
Consider the problem \eqref{nls} where the cubic nonlinearity $C$ satisfies the 
assumptions (H1) and (H2). Assume that the initial data $u_0$ is small, 
\[
\|u_0\|_{L^2_x} \leq \epsilon \ll 1,
\] 
Then the solution $u$ exists on on a time interval $I_\epsilon = [0,c \epsilon^{-8}]$ and 
has the following properties for every interval $I \subset I_\epsilon$ of size $|I| \leq \epsilon^{-6}$:

\begin{enumerate}[label=(\roman*)]
\item Uniform $L^2$ bound:
\begin{equation}\label{main-L2}
\| u \|_{L^\infty_t(I_\epsilon; L^2_x)} \lesssim \epsilon.
\end{equation}

\item Strichartz bound:
\begin{equation}\label{main-Str}
\| u \|_{L^6_{t,x}(I \times \R)} \lesssim \epsilon^\frac23.
\end{equation}

\item Bilinear Strichartz bound:
\begin{equation}\label{main-bi}
\| \partial_x (u \bar u(\cdot+x_0))\|_{L^2_t (I;H_x^{-\frac12})} \lesssim \epsilon^2,
\qquad x_0 \in \R.
\end{equation}
\end{enumerate}
\end{theorem}
The local well-posedness in $L^2$ for  the problem \eqref{nls} 
was already proved in \cite{IT-noloc}, so the emphasis here and later in the proof is on the lifespan bound and the long time estimates in the theorem.

We remark that the intermediate time-scale $\epsilon^{-6}$ does not have an intrinsic meaning from a scaling perspective, but is instead connected to the unit frequency-scale 
which is implicit in (H1), and which motivates the frequency decomposition on the 
unit frequency-scale which is used in the proof of both the present result and the earlier result in \cite{IT-noloc}. One could also use a smaller frequency-scale for this decomposition, which in turn corresponds to a smaller size for $|I|$. This is however not needed in the proof of the 
$\epsilon^{-8}$ result, so, in order to avoid cluttering the theorem, we omit the details.
But the interested reader should see Remark~\ref{r:scales} below.

A natural question to ask is whether this result is optimal.
On one hand, the bounds \eqref{main-Str} and \eqref{main-bi} are sharp for the cubic NLS,
and likely for any focusing flow (i.e. which satisfies (H4) at least in some region); this is discussed in Section~\ref{s:sharp}.  But the $\epsilon^{-8}$ lifespan bound is not optimal for any flow satisfying our hypotheses. Indeed, if for instance the $L^2$ norm is conserved (as is the case for
the focusing NLS) then global well-posedness follows.  However, we conjecture that
\begin{conjecture}
The  result in Theorem~\ref{t:focusing} is sharp for generic focusing flows satisfying our hypotheses.
\end{conjecture}

It it also interesting to see how our theorem applies to the focusing cubic NLS problem. There we can also consider large data simply by scaling.
Global well-posedness in $L^2$ is relatively straightforward there, but some of the estimates we prove are new:

\begin{theorem}
Consider the focusing 1-d cubic NLS problem \eqref{nls3}(-) with $L^2$ initial data $u_0$.
Then the global solution $u$ satisfies the following bounds in all intervals $I$ 
with size $|I| > \|u_0\|_{L^2_x}^{-4}$:

\begin{enumerate}[label=(\roman*)]
\item Uniform $L^2$ bound:
\begin{equation}\label{main-L2-model}
\| u \|_{L^\infty_t L^2_x} \lesssim \|u_0\|_{L^2_x}.
\end{equation}

\item Strichartz bound:
\begin{equation}\label{main-Str-model}
\| u \|_{L^6_{t,x}(I\times \R)} \lesssim  \|u_0\|_{L^2_x} ( I \|u_0\|_{L^2_x}^4)^\frac16.
\end{equation}

\item Bilinear Strichartz bound:
\begin{equation}\label{main-bi-model}
\| \partial_x |u|^2\|_{L^2_t (I;\dot H_x ^{-\frac12} + c L^2_x) } \lesssim  \|u_0\|_{L^2_x}^2, \qquad c^{2} = \| u_0\|_{L^2_x}^{2} (|I|   \|u_0\|_{L^2_x}^4).
\end{equation}
\end{enumerate}
\end{theorem}
Here the Strichartz estimates are fairly easy to obtain directly, but the bilinear Strichartz bounds are new. Returning to our discussion after Theorem~\ref{t:focusing}, the numerology 
in this application helps clarify the earlier comment about the choice of the time-scales
in Theorem~\ref{t:focusing}:
\begin{remark}\label{r:scales}
In the context of Theorem~\ref{t:focusing}, a scaling argument shows that the bounds  
\eqref{main-Str-model} and \eqref{main-bi-model} hold for all intervals $I$
so that $\epsilon^{-4} \leq |I| \leq \epsilon^{-6}$. 
\end{remark}

One may gain further insights into our result for focusing problems by 
comparing it with our earlier result in \cite{IT-noloc}.

\begin{theorem}[\cite{IT-noloc}]\label{t:defocusing}
Under the above assumptions (H1), (H2) and (H3) on the symbol of the cubic form $C$, small initial data 
\[
\|u_0\|_{L^2_x} \leq \epsilon \ll 1,
\] 
yields a unique global solution $u$ for \eqref{nls}, which satisfies the following bounds:

\begin{enumerate}[label=(\roman*)]
\item Uniform $L^2$ bound:
\begin{equation}\label{main-L2-de}
\| u \|_{L^\infty_t L^2_x} \lesssim \epsilon.
\end{equation}

\item Strichartz bound:
\begin{equation}\label{main-Str-de}
\| u \|_{L^6_{t,x}} \lesssim \epsilon^\frac23.
\end{equation}

\item Bilinear Strichartz bound:
\begin{equation}\label{main-bi-de}
\| \partial_x (u \bar u(\cdot+x_0))\|_{L^2_t H_x^{-\frac12}} \lesssim \epsilon^2,
\qquad x_0 \in \R.
\end{equation}
\end{enumerate}
\end{theorem}

One may observe here that the estimates are similar in the two cases, and the only difference
is the time-scale on which the estimates hold: in the defocusing case this is global, while in the focusing case it is finite and depends on the solution size. For this reason, the 
proofs of Theorem~\ref{t:defocusing} and \eqref{t:focusing} are closely related, and we will take advantage of this within the proof.

For the convenience of the reader, we also recall the main ideas in the proof of the last theorem in \cite{IT-noloc}, which are 
equally employed here:

\medskip

\emph{ 1. Energy estimates via density flux identities.} This is a classical idea in pde's, and particularly in the study of conservation laws. The novelty in \cite{IT-noloc}  is that this analysis is carried out in  a nonlocal setting, where both the densities and the fluxes involve translation invariant multilinear forms.
The densities and the fluxes are not uniquely determined
here, so careful choices need to be made.
\medskip

\emph{2. The use of energy corrections.}
This is an idea originally developed in the context of the so called 
I-method~\cite{I-method} or more precisely the second generation I-method \cite{I-method2}, whose aim was to construct more accurate almost conserved quantities. In \cite{IT-noloc} this idea is instead implemented at the level of density-flux
identities. 
\medskip

\emph{3. Interaction Morawetz bounds.} These were originally developed in the context of the three-dimensional NLS problems by Colliander-Keel-Stafillani-Takaoka-Tao in \cite{MR2053757},
and have played a fundamental role in the study of many nonlinear Schr\"odinger flows, see e,g. \cite{MR2415387,MR2288737}, and also for one-dimensional quintic flows in the work of Dodson~\cite{MR3483476,MR3625190}. Our take on this is somewhat closer to the one-dimensional  approach of Planchon-Vega \cite{PV},
though recast in the setting and language of nonlocal multilinear 
forms.
\medskip

\emph{4. Tao's frequency envelope method.} This is used as 
a way to accurately track the evolution of the energy distribution across frequencies. Unlike the classical implementation relative to dyadic Littlewood-Paley decompositions, in \cite{IT-noloc} we adapt and refine this notion for lattice decompositions instead. 
This is also very convenient as a bootstrap tool,  see e.g. Tao~\cite{Tao-WM}, \cite{Tao-BO} but 
with the added twist of also bootstrapping bilinear Strichartz bounds, as in the authors' paper \cite{IT-BO}.

\subsection{ An outline of the paper} 
To a large extent, the proof of our main result mirrors the proof of the global result in the defocusing case in \cite{IT-noloc}.  The primary difference is in how the $L^6$ Strichartz norms are handled, both globally
and in a frequency localized setting. 

Section~\ref{s:Morawetz} reviews two of the main ideas in \cite{IT-noloc}, namely the construction of modified 
density-flux identities for the mass/momentum in a frequency localized setting, as well as the interaction Morawetz identities associated to 
those density-flux relations.

The proof of the main $L^6$ Strichartz 
bounds and the bilinear $L^2$ estimates
is done in a frequency localized setting, using a bootstrap argument  based on a frequency decomposition on the unit scale, where the components are measured using a maximal frequency envelope, another notion introduced in \cite{IT-noloc}. This is described in Section~\ref{s:boot}, and leads 
to dispersive bounds on the $\epsilon^{-6}$ time-scale.
It is within this argument where the $L^6$ norms
are treated differently from the defocusing case.

In order to advance from the $\epsilon^{-6}$ to the $\epsilon^{-8}$ time-scale it suffices to propagate 
the  $L^2$ bound (i.e. the mass) on the larger time-scale. However the mass is not a conserved quantity,
so instead it is better to propagate the bounds  for the modified mass. This analysis is carried out in Section~\ref{s:energy}.

Finally, in the last section of the paper we discuss the optimality of our result, or rather the optimality 
of the $L^6$ and the bilinear $L^2$ estimates on the 
$\epsilon^{-6}$ time-scale. This is done by considering the obvious enemies, namely the solitons, in the focusing NLS context.

\subsection{Acknowledgements} 
 The first author was supported by the Sloan Foundation, and by an NSF CAREER grant DMS-1845037. The second author was supported by the NSF grant DMS-2054975 as well as by a Simons Investigator grant from the Simons Foundation.

\

\section{Density-flux and interaction Morawetz identities}
\label{s:Morawetz}

A key role in the proof of the results in both \cite{IT-noloc}
and in the present paper is played by the approximate conservation laws for the mass and the momentum. Rather than considering them directly, we instead consider several improvements:
\begin{itemize}
\item the conservation laws are written in density-flux form,
rather than as integral identities, where both the densities 
and the fluxes are multilinear forms.

\item we improve the accuracy of these conservation laws 
by using well chosen quartic corrections for both the densities
and the fluxes, with $6$-linear errors.

\item we use these densities and fluxes not only globally in frequency, but also in a frequency localized setting.
\end{itemize}
The aim of this section is to provide an overview  of these density-flux identities, following the set-up of \cite{IT-noloc}. We conclude the section with an overview 
of the interaction Morawetz identities obtained in \cite{IT-noloc} from the above density flux identities.

\subsection{Resonances and multilinear forms}

A key role in our analysis is played by four wave resonances.
Given three input frequencies $(\xi_1,\xi_2,\xi_3)$ in the cubic nonlinearity $C$, the output is at frequency 
\[
\xi_4 = \xi_1 - \xi_2 + \xi_3.
\]
The three wave interaction is resonant if 
\[
\xi_4^2 = \xi_1^2 - \xi_2^2 + \xi_3^2.
\]
To rewrite these relations in a symmetric fashion we use the notations
\[
\Delta^4 \xi = \xi_1 - \xi_2 + \xi_3-\xi_4,
\qquad \Delta^4 \xi^2 = \xi_1^2 - \xi_2^2 + \xi_3^2 -\xi_4^2,
\]
The first expression is Galilean invariant but not the second, which is why we also use the adjusted, Galilean invariant expression
\[
\tDelta^4 \xi^2 = \Delta^4 \xi^2 - 2 \xi_{avg} \Delta^4 \xi,
\]
where $\xi_{avg}$ represents the average of the four frequencies.

With these notations, the resonant set is described as 
\[
\calR = \{ [\xi] = (\xi_1,\xi_2,\xi_3,\xi_4); \ \Delta^4 \xi = 0,
\Delta^4 \xi^2 = 0 \}
\]
which can be explicitly characterized as 
\[
\calR = \{ [\xi]; (\xi_1,\xi_3) = (\xi_2,\xi_4) \}
\]
Quadruples in the resonant set can be described by two parameters, 
namely
\begin{equation}\label{xihm}
\begin{aligned}
\xih = & \ \max  \{|\xi_1 - \xi_{2}|+|\xi_3-\xi_4|, |\xi_1 - \xi_{4}|+|\xi_3-\xi_2| \},
\\
\xim =& \  \min  \{|\xi_1 - \xi_{2}|+|\xi_3-\xi_4|, |\xi_1 - \xi_{4}|+|\xi_3-\xi_2| \},
\end{aligned}
\end{equation}
These distance parameters are carefully defined so that they can also be used outside the resonant set to characterize frequency quadruples. This is very useful in the density-flux relations
later on.

\subsection{The conservation of mass}

The starting point of the analysis in \cite{IT-noloc}
is to consider energy estimates for our flow from a density-flux perspective.
In the simplest case, we start with the mass density
\[
M(u,\bar u) = |u|^2,
\]
whose linear flux is given by the momentum
\[
P(u,\bar u) = 2i \Im (u \partial_x \bar u). 
\]
These can be viewed as translation invariant bilinear forms 
with symbols 
\[
m(\xi,\eta) = 1, \qquad p(\xi,\eta) = \xi+\eta.
\]
Integrating the densities we obtain the familiar mass 
and momentum,
\[
\bM(u) = \int_{\mathbb{R}} M(u,\bar u) \, dx, \quad \bP(u) = \int_{\R} M(u,\bar u) \, dx.
\]

At the nonlinear level, we have the density flux relation
\[
\partial_t M(u,\bar u) = \partial_x P(u,\bar u) + C^4_m(u,\bar u,u,\bar u),
\]
where $C^4_m$ is a symmetric translation invariant real multilinear form which depends on our cubic nonlinearity 
$C$. A key observation in \cite{IT-noloc} is that, under the conservative assumption (H2) on the nonlinearity, the mass density admits a quartic correction which is accurate to sixth order.
Precisely, the correction has the form
\begin{equation}\label{m-sharp}
\ms(u) = M(u) + B^4_m(u,\bar u, u,\bar u),
\end{equation}
and the associated density flux relation has the form
\begin{equation}\label{dens-flux-m}
\partial_t \ms(u) = \partial_x (P(u) + R^4_m(u,\bar u, u,\bar u)) +
R^6_m(u,\bar u,  u,\bar u, u,\bar u) .
\end{equation}
with suitable translation invariant multilinear forms
$R^4_m$ and $R^6_m$. The corresponding integral corrected mass is
\[
\bM^\sharp(u) = \int_{\R} \ms(u) \, dx.
\]

The choice of the symbols $b^4_m$ and $r^4_m$ above
depends on the behavior of $c^4_m$ near the resonant set $\calR$,
precisely they have to solve the division problem 
\begin{equation}\label{choose-R4m}
c^4_m +  i  \Delta^4 \xi^2 \,  b^4_m = i \Delta^4 \xi\, r^4_m .
\end{equation}
This is possible due to our condition (H2), which implies 
that $c^4_m = 0$ on $\calR$.
But the choice is not uniquely determined, so it is important to make a good one, i.e. which insures good symbol bounds. To achieve
this, in \cite{IT-noloc} we
decompose the phase space for frequency quadruples into three overlapping regions which can be separated using cutoff functions which are smooth on the unit scale:
\begin{enumerate}[label=\roman*)]
\item The full division region,
\[
\Omega_1 = \{ \xim \lesssim 1 \},
\]
which represents a full unit size neighbourhood of the resonant set $\calR$.
\item The region 
\[
\Omega_2  = \{ 1+ |\Delta^4 \xi| \ll \xim \}, 
\]
where $\tDelta^4 \xi^2$  must be elliptic, $\tDelta^4 \xi^2 \approx \xih \xim$,
and thus we will favor division by the symbol $\tDelta^4 \xi^2$.

\item The region
\[
\Omega_3 =  \{ 1 \ll \xim \lesssim |\Delta^4 \xi|  \},
\]
we will instead divide by $\Delta^4 \xi$; this 
is compensated by the relatively small size of this region.
\end{enumerate}

Since we will also  need this in the present paper, we state the result in the following 

\begin{proposition}[\cite{IT-noloc}]
Assume that the nonlinearity $C$ satisfies the conditions (H1), (H2). Then there exist 
multilinear forms $B^4_m$ $R^4_m$ and  $R^6_m$ so that the  relation  \eqref{dens-flux-m} holds for solutions $u$ to \eqref{nls}, and so that the symbols $c^4_m$ and $r^4_m$ satisfy the bounds
\begin{enumerate}[label = \roman*)]
    \item Size 
\begin{equation}
\begin{aligned}
|\partial^\alpha  r^4_m|  \lesssim & \  \frac{1}{\la \xim \ra} , 
\\
|\partial^\alpha b^4_m|  \lesssim & \  \frac{1}{\la\xih\ra \la\xim\ra}.
\end{aligned}   
\end{equation}    
\item Support:
$b^4$ is supported in $\Omega_1 \cup \Omega_2$ and $\tilde r^4$ is supported in $\Omega_1 \cup \Omega_3$.
\end{enumerate}
In addition, we have the fixed time bound
\begin{equation} \label{small-correction}
\left| \int_{\R} B^4_m(u,u,u,u) \, dx \right| \lesssim \|u\|_{L^2_x}^4. 
\end{equation}
\end{proposition}
We make several remarks concerning this result:
\begin{itemize}
    \item No bound for the symbol $r^6_m$ is provided in the proposition. This is because $r^6_m$ is obtained directly as the contribution of $C$ to the time derivative of $b^4_m$.
\item This proposition is a consequence of Lemma 4.1 and Lemma 7.1 in \cite{IT-noloc}.
\item 
The defocusing hypothesis (H3), which is used for the final result 
in \cite{IT-noloc}, plays no role here.

\item The estimate \eqref{small-correction} shows that the mass correction is perturbative for as long as the solution $u$ remains small in $L^2$.

\item A similar analysis applies for the momentum conservation law. But the counterpart of the above Proposition for the momentum is less useful directly, and instead it is used in \cite{IT-noloc} only in a frequency localized context. 
\end{itemize}

\subsection{Frequency localized density-flux identities}

Instead of relying on the more standard Littlewood-Paley decomposition, the analysis in \cite{IT-noloc}
uses a frequency decomposition on the unit scale in frequency. Given any integer $j$,
we will use localized versions of the mass in a unit size region around $j$. More generally,
for an interval $A \subset \Z$, we use a symbol $a_0$ which is frequency localized in a unit neighbourhood of $A$. At the level of bilinear forms, we will use the symbol 
\[
a(\xi,\eta) = a_0(\xi) a_0(\eta).
\]
Corresponding to such $a$ we define  quadratic 
localized mass, momentum and energy densities by
the symbols
\[
m_a(\xi,\eta) = a(\xi,\eta), \qquad p_a(\xi,\eta) = (\xi+\eta)
a(\xi,\eta), \qquad e_a(\xi,\eta) = (\xi+\eta)^2
a(\xi,\eta).
\]
The associated bilinear forms are denoted by $M_a$, $P_a$, respectively $E_a$.
If $A = \{ j\}$ then we simply replace the subscript $a$ with $j$.

\medskip

It is shown in \cite{IT-noloc}, again under the assumptions (H1) and (H2), 
that one may find quartic corrections $\ms_a$ and $\ps_a$
of the form

\begin{equation}\label{ma-sharp}
\ms_a(u) = M_a(u) + B^4_{m,a}(u,\bar u, u,\bar u),
\end{equation}
\begin{equation}\label{pa-sharp}
\ps_a(u) = P_a(u) + B^4_{p,a}(u,\bar u, u,\bar u),
\end{equation}
for which we obtain density-flux identities akin to
\eqref{dens-flux-m},  namely 
\begin{equation}\label{dens-flux-ma}
 \partial_t \ms_a(u) = \partial_x(P_{a}(u)
+ R^4_{m,a}(u)) +  R^6_{m,a}(u),
\end{equation}
and 
\begin{equation}\label{dens-flux-pa}
 \partial_t \ps_a(u) = \partial_x(E_{a}(u)
+ R^4_{p,a}(u)) +  R^6_{p,a}(u).
\end{equation}

We will consider these relations together with their Galilean shifts
obtaining relations of the form 
\begin{equation}\label{dens-flux-maG}
(\partial_t - 2 \xi_0 \partial_x)\ms_a(u) = \partial_x(P_{a,\xi_0}(u)
+ R^4_{m,a,\xi_0}(u)) +  R^6_{m,a,\xi_0}(u),
\end{equation}
respectively 
\begin{equation}\label{dens-flux-paG}
(\partial_t - 2 \xi_0 \partial_x)\ps_{a,\xi_0} (u) = \partial_x(E_{a,\xi_0}(u)
+ R^4_{p,a,\xi_0}(u)) +  R^6_{p,a,\xi_0}(u).
\end{equation}

These correspond to the algebraic division relations
\begin{equation}\label{choose-R4m-a}
c^4_{m,a} +  i  \Delta^4 (\xi-\xi_0)^2 b^4_{m,a} = i \Delta^4 \xi r^4_{m,a,\xi_0} ,
\end{equation}
respectively
\begin{equation}\label{choose-R4p-a}
c^4_{p,a,\xi_0} +  i  \Delta^4 (\xi-\xi_0)^2 b^4_{p,a,\xi_0} = 
i \Delta^4\xi r^4_{p,a,\xi_0},
\end{equation}
where $c^4_{m,a}$ and $c^4_{p,a,\xi_0}$ are the density-flux sources
corresponding to uncorrected mass, respectively momentum.

The symbols above are connected in the obvious way. Precisely,
we have 
\begin{equation}
  r^4_{m,a,\xi_0} =   r^4_{m,a} - 2\xi_0 b^4_{m,a},
 \qquad 
\end{equation}
and
\begin{equation}
\ps_{a,\xi_0} = \ps_a - 2\xi_0 \ms_a, \qquad b^4_{p,a,\xi_0}
= b^4_{p,a} - 2 \xi_0 b^4_{m,a}, 
\end{equation}
and finally 
\begin{equation}
r^4_{p,a,\xi_0} = r^4_{p,a} - 2\xi_0 b^4_{p,a} - 2 \xi_0 r^4_{m,a,\xi_0}.
\end{equation}

To use these density flux relations we need to have appropriate 
bounds for our symbols:

\begin{proposition} 
\label{p:symbols}
Let $J \subset \R$ be an interval of length $r$, and $d(\xi_0,J) \lesssim r$.  Assume that $a$ is supported in $J \times J$, with bounded and uniformly smooth symbol. Then the relations \eqref{dens-flux-maG}
and \eqref{dens-flux-paG} hold with symbols $b^{4}_{m,a}$, $b^{4}_{p,a,\xi_0}$, $r^4_{m,a,\xi_0}$ and $r^4_{p,a,\xi_0}$
which can be chosen to have the following properties:

\begin{enumerate}[label=\roman*)]
    \item Support: they are all supported in the region where 
    at least one of the frequencies is in $J$.
    \item Size: 
\begin{equation}
 |b^{4}_{m,a}| \lesssim \frac{1}{\la \xih \ra \la \xim \ra }, 
 \qquad
 |b^{4}_{p,a,\xi_0}| \lesssim \frac{r}{\la \xih \ra \la \xim \ra },
\end{equation}    
\begin{equation}
\begin{aligned}
 |r^{4}_{m,a,\xi_0}| \lesssim  & \ \frac{1}{\la \xim  \ra } 1_{\Omega_1 \cup \Omega_2}
 + \frac{r}{\la \xih \ra \la \xim \ra }
   1_{\Omega_1 \cup \Omega_3},
 \\
 |R^{4}_{p,a,\xi_0}| \lesssim & \ \frac{r}{\la \xim  \ra } 1_{\Omega_1 \cup \Omega_2}
 + \frac{r^2}{\la \xih \ra \la \xim \ra }
   1_{\Omega_1 \cup \Omega_3}.
 \end{aligned}  
\end{equation}     
  \item Regularity:  similar bounds hold for all derivatives.
\end{enumerate}
\end{proposition}

 This is Proposition~4.3 in \cite{IT-noloc}. 

\subsection{ Interaction Morawetz identities}
One way the desity flux relations above are used is to obtain
more accurate bounds for the mass propagation. However, 
another way to use them is via interaction Morawetz identities,
which yield bilinear estimates for the interaction of 
different frequency portions of the solutions, or even the self-interaction of unit frequency portions of solutions.

Given two frequency intervals $A$ and $B$ and corresponding
mass/momentum modified density associated to these intervals
for two solutions $u,v$ to \eqref{nls},
in \cite{IT-noloc} we define the associated \emph{interaction Morawetz 
functional} by
\begin{equation}\label{interaction-bi}
 \bI_{AB} = \iint_{x > y} \ms_a(u)(x) \ps_{b,\xi_0}(v)(y)    - \ps_{a,\xi_0}(u)(x) \ms_b(v)(y)\, dx dy ,
\end{equation}
where we add several remarks:
\begin{itemize}
\item In applications, for the second solution $v$ we will simply choose  a spatial translation of the solution $u$.
\item  A velocity parameter $\xi_0$ is introduced on the right, 
but the functional does not depend on $\xi_0$. This parameter 
plays a role, however, in estimating the time derivative of 
$\bI_{AB}$, and will be chosen to be close to the sets $A$ and $B$.
\item Such interaction functionals will be used in two settings:
\begin{enumerate}[label=(\roman*)]
\item  the separated case  $ 1 \ll |A| \approx |B| \approx dist(A,B)$, and 
\item  the self-interaction case $|A| + |B|+  dist(A,B) \lesssim 1$.
\end{enumerate}
\end{itemize}

The time derivative of the interaction Morawetz 
functional is computed using the frequency localized mass density-flux \eqref{dens-flux-maG}  and the corresponding momentum density-flux \eqref{dens-flux-paG}. 
This yields a localized \emph{interaction Morawetz identity},
\begin{equation}\label{interaction-xi-AB}
\frac{d}{dt} \bI_{AB} =  \bJ^4_{AB} + \bJ^6_{AB} + \bJ^8_{AB} + \bK^8_{AB},
\end{equation}
where the terms on the right are described as follows:
\begin{enumerate}[label=\alph*)]
\item the quartic contribution $\bJ^4_{a}$ is 
\[
\bJ^4_{AB} = \int_{\R} M_a(u)(x) E_{b,\xi_0}(v)(x) + M_b(v)(x) E_{a,\xi_0}(u)(x)-  2 P_{a,\xi_0}(u)(x) P_{b,\xi_0}(v)(x)\, dx,
\]
and is used to capture bilinear $L^2$ bounds.

\item The sixth order term $\bJ^6_{a}$ has the form
\[
\bJ^6_{AB} =   \int_{\R} -( P_{a,\xi_0} B^4_{p,b,\xi_0}
+ P_{b,\xi_0} R^4_{m,a,\xi_0}) + ( 
M_a R^4_{p,b,\xi_0} +E_{b,\xi_0} B^4_{m,a,\xi_0})      - \text{symmetric} \, dx,
\]
where in the symmetric part we interchange both the indices $a,b$ and the functions $u,v$. In the defocusing case its symbol has a
favorable sign on the diagonal, and is used to capture the $L^6$ bound in the self-interaction case. Here, it is estimated perturbatively.

\item The eight-linear term
\[
\bJ^8_{AB} =   \int_{\R} - R^4_{m,a,\xi_0} B^4_{p,b,\xi_0}
 +  
 B^4_{m,a,\xi_0} R^4_{p,b,\xi_0}  - \text{symmetric} \,      dx .
\]
plays a perturbative role.

\item The $8^+$-linear term $\bK^8_{AB}$ has the form
\[
\bK^8_{AB} = \iint_{x > y} \ms_a(x)  R^6_{p,a,\xi_0} + \ps_{a,\xi_0} R^6_{m,a,\xi_0} - \text{symmetric} \, dx dy.
\]
This is a double integral, which also includes  a $10$-linear term.
It is also estimated perturbatively.
\end{enumerate}

\section{Strichartz and bilinear \texorpdfstring{$L^2$}{} bounds}
\label{s:boot}

To obtain estimates for $L^2$ solutions $u$ to \eqref{nls},
a unit scale frequency decomposition is needed,
\[
u = \sum_{k \in \Z} u_k, \qquad u_k := P_k u,
\]
where $P_k$ are multipliers with smooth symbols 
localized in a unit neighbourhood of the integer frequency $k$.
To measure the components $u_k$ we use a frequency envelope $\{c_k\} \in \ell^2$ in order to transfer bounds from the initial data to the solutions.

These frequency envelopes are chosen to satisfy an adapted version 
of the slowly varying property, originally introduced by Tao \cite{Tao-WM} in the context of dyadic decompositions. Such a property is needed in order
to account for the nonlinear leakage of energy between nearby frequencies.

We recall the frequency envelope set-up in \cite{IT-noloc}, associated to lattice decompositions:

\begin{definition}
A lattice frequency envelope $\{c_k\}$ is said to  
have the maximal property if 
\begin{equation}
Mc \leq C c    ,
\end{equation}
where $Mc$ represents the maximal function of $c$.
Here $C$ is a universal constant.
\end{definition}

Frequency envelopes that have 
this property will be called \emph{admissible}.
The proof of the $L^6$ Strichartz and the bilinear $L^2$ 
bounds will be phrased as a bootstrap argument relative to 
an admissible  frequency envelope for the initial data.

\begin{theorem}\label{t:boot}
 Let $u \in C[0,T;L^2]$ be a solution for the equation \eqref{nls} with initial data $u_0$  which has $L^2$ size at most $\epsilon$. 
Let $\{c_k\}$ be a maximal frequency envelope for the initial data  in $L^2$, also of size $\epsilon$,
\[
\| u_{0k}\|_{L^2} \lesssim \epsilon c_k.
\]
Assume that 
 \begin{equation} \label{T-bound}
    T \ll \epsilon^{-6}. 
 \end{equation}
Then the solution $u$ satisfies  the following bounds in $[0,T]$:
\begin{enumerate}[label=(\roman*)]
\item Uniform frequency envelope bound:
\begin{equation}\label{uk-ee}
\| u_k \|_{L^\infty_t L^2_x} \lesssim \epsilon c_k,
\end{equation}
\item Localized Strichartz bound:
\begin{equation}\label{uk-se}
\| u_k \|_{L_{t,x}^6} \lesssim (\epsilon c_k)^\frac23,
\end{equation}
\item Localized Interaction Morawetz:
\begin{equation}\label{uk-bi}
\| \partial_x |u_k|^2  \|_{L^2_{t,x}} \lesssim \epsilon^2 c_k^2,
\end{equation}
\item Transversal bilinear $L^2$ bound:
\begin{equation} \label{uab-bi}
\| \partial_x(u_{A} \bu_B(\cdot+x_0))  \|_{L^2_{t,x}} \lesssim \epsilon^2 c_A c_B\, \la dist(A,B)\ra^{\frac12},
\end{equation}
 for all $x_0 \in \R$ whenever $|A| + |B| \lesssim  \la dist(A,B)\ra $.
\end{enumerate}
\end{theorem}
This proposition mirrors a similar result in Section 7 of \cite{IT-noloc},
with two key differences. On one hand we drop the defocusing assumption (H3), and on the other hand we limit the size of the time interval in \eqref{T-bound}.

The $L^6$ Strichartz estimates \eqref{main-Str} and the bilinear $L^2$ bounds in \eqref{main-bi} follow from 
the estimates in the 
above proposition, by the same arguments as those in Section~8 of \cite{IT-noloc}.

To prove this theorem, we make a bootstrap assumption where we assume the same bounds but with a worse constant $C$, as follows:

\begin{enumerate}[label=(\roman*)]
\item Uniform frequency envelope bound,
\begin{equation}\label{uk-ee-boot}
\| u_k \|_{L^\infty_t L^2_x} \lesssim C \epsilon c_k,
\end{equation}
\item Localized Strichartz bound,
\begin{equation}\label{uk-se-boot}
\| u_k \|_{L^6_{t,x}} \lesssim C(\epsilon c_k)^\frac23,
\end{equation}
\item Localized Interaction Morawetz,
\begin{equation}\label{uk-bi-boot}
\| \partial_x |u_k|^2  \|_{L^2_{t,x}} \lesssim  C\epsilon^2 c_k^2,
\end{equation}
\item Transversal Interaction Morawetz, 
\begin{equation} \label{uab-bi-boot}
\| \partial_x (u_{k_1} \bar u_{k_2}(\cdot+x_0))  \|_{L^2_{t,x}} \lesssim C \epsilon^2 c_{k_1} c_{k_2} \la k_1-k_2\ra^{\frac12}
\end{equation}
uniformly for all $x_0 \in \R$.
\end{enumerate}

Then we seek to improve the constant in these bounds. The gain will come from the fact that the $C$'s will always come paired
either with extra $\epsilon$ factors, or with $T\epsilon^6$ 
factors.

To a large extent the proof largely repeats the proof of the corresponding result in \cite{IT-noloc}, so we review the steps
and expand the portion where the argument differs here.

\bigskip

\emph{ STEP 1:} The proof of the energy bound \eqref{uk-ee}.
This is done by integrating the density flux relation 
\eqref{dens-flux-ma} for the localized mass $\bM_k(u)$. The argument 
in \cite{IT-noloc} applies unchanged.

\bigskip

\emph{ STEP 2:} The proof of the $L^6$ Strichartz bound \eqref{uk-ee}. In \cite{IT-noloc} this is proved together with 
\eqref{uk-bi} by integrating the interaction Morawetz identity applied to the pair $(u_k,u_k)$. However, in the focusing case
the sign of the $u_k^6$ contribution changes, and the same argument no longer applies. 

Instead, here we will estimate the $L^6$ norm directly using an interpolation argument. Precisely, for the function $v_k = |u_k|^2$, using \eqref{uk-ee-boot} and \eqref{uk-bi-boot}, we have the bounds
\[
\| v_k\|_{L^\infty_t L^1_x} \lesssim C^2 \epsilon^2 c_k^2 ,
\]
respectively
\[
\| v_k\|_{L^2_t \dot H^1_x} \lesssim C^2 \epsilon^2 c_k^2 .
\]
We interpolate between the two estimates in homogeneous Sobolev spaces, with weights $5/9$ 
and $4/9$. We obtain
\[
\| v_k\|_{L^{\frac92}_t \dot{W}^{\frac{4}{9}, \frac{9}{7}}_x}\lesssim \Vert v_k\Vert_{L^\infty_t L^1_x}^{\frac{5}{9}}\|  v_k\|_{L^2_t \dot{H}^1_x}^{\frac{4}{9}} \lesssim C^2 \epsilon^2 c_k^2.
\]
By Sobolev embeddings  $\dot{W}^{\frac{4}{9}, \frac{9}{7}}$ embeds in $L^3$ so, using also H\"older's inequality with respect to time, we obtain
\[
\Vert v_k\Vert_{L^3_t L_x^3}\lesssim T^{\frac{1}{9}} \Vert v_k\Vert _{L^{\frac{9}{2}}L^3} \lesssim T^{\frac{1}{9}}\| v_k \|_{L^{\frac{9}{2}} \dot{W}^{\frac{4}{9},\frac97}}\lesssim C^2 T^{\frac{1}{9}}\epsilon^2 c_k^2 = C^2 (T \epsilon^6)^\frac19 
\epsilon^{\frac43} c_k^2.
\]
This implies the desired Strichartz bound \eqref{uk-se} under the
time constraint \eqref{T-bound}.

\bigskip

\emph{ STEP 3:} The proof of the bilinear $L^2$ bound \eqref{uab-bi}. This is again exactly as in \cite{IT-noloc},
by applying the interaction Morawetz identity to the functions
$(u_A, u_B(\cdot+x_0))$.  We note that the $L^6$ bound is used as an input in this proof, and the defocusing assumption (H2) is not needed. Here we view \eqref{uk-bi} as a special case of \eqref{uab-bi}, and no longer in conjunction with   \eqref{uk-ee}.

\section{Long time energy estimates}
\label{s:energy}
The frequency envelope bounds in Theorem~\ref{t:boot}
provide us with uniform energy bounds on the $\epsilon^{-6}$
time-scale, and so they do not suffice in order to prove our main result in Theorem~\ref{t:focusing}, which is on the $\epsilon^{-8}$
time-scale. To fill in this gap, we will prove a direct energy estimate on the $\epsilon^{-8}$ time-scale. Precisely, we will show 
the following:

\begin{proposition}
Let $u$ be an $L^2$ solution for \eqref{nls} in a time interval $[0,T]$. Assume that the initial data for \eqref{nls} satisfies
\begin{equation}
\|u_0\|_{L_x^2} \leq \epsilon  \ll 1 , 
\end{equation}
and that $T \ll \epsilon^{-8}$. 
Then the solution $u$ satisfies
\begin{equation}
\| u\|_{L_t^\infty (0,T;L_x^2)} \leq 4 \epsilon.
\end{equation}
\end{proposition}

Once we have this proposition, a continuity argument based on the local well-posedness for \eqref{nls} in $L^2$ implies 
Theorem~\ref{t:focusing}.

\begin{proof}
It suffices to prove that the conclusion holds assuming that we have the bootstrap assumption
\begin{equation}
\| u\|_{L^\infty (0,T;L^2)} \leq 8 \epsilon.
\end{equation}

Instead of tracking directly the mass $\bM(u) = \|u\|_{L_x^2}^2$,
it is more efficient to work with the modified mass 
\[
\bM^\sharp(u) = \int_\R \ms(u) \, dx.
\]
In view of the bound \eqref{small-correction}, we have 
\[
\bM(u) = \bM^\sharp(u) + O(\epsilon^4).
\]
Since $\epsilon \ll 1$, we have 
\begin{equation}\label{bm-0}
\bM(u)(0) \leq 2 \epsilon^2,
\end{equation}
and it suffices to show  that 
\begin{equation}\label{bm-t}
  \bM(u)(t) \leq 4 \epsilon^2, \qquad t \in [0,T]  .
\end{equation}
In view of Proposition~\ref{dens-flux-m}, the time evolution 
of $\bM^\sharp$ is given by
\begin{equation}
\frac{d}{dt}    \bM^\sharp(u) = \int_{\R} R^6_m(u,\bu,u,\bu,u,\bu)\, dx. 
\end{equation}

To bound its growth, we use the following
\begin{lemma}
\label{l:R6-AB}
Assume that the bounds \eqref{uk-ee}-\eqref{uab-bi} for $u$ 
hold in a time interval $[0,T]$. Then we have 
\begin{equation}\label{R6-m-bd}
\|R^6_{m}\|_{L^1_{t,x}([0,T]\times \R)} \lesssim \epsilon^4.
\end{equation}
\end{lemma}
This is Lemma 7.3 in \cite{IT-noloc}.
We apply this lemma on time intervals of size $\epsilon^{-6}$,
where the bounds \eqref{uk-ee}-\eqref{uab-bi} hold in view of Theorem~\ref{t:boot}, and then add up the results. Then for $T > \epsilon^{-6}$ we get
\begin{equation}\label{R6-m-bd-a}
\|R^6_{m}\|_{L^1_{t,x}([0,T]\times \R)} \lesssim \epsilon^4 (T \epsilon^{6}) .
\end{equation}
Hence for $T \ll \epsilon^{-8}$ we arrive at
\begin{equation}\label{R6-m-bd-b}
\|R^6_{m}\|_{L^1_{t,x}([0,T]\times \R)} \ll \epsilon^2 .
\end{equation}
This allows us to obtain \eqref{bm-t} from \eqref{bm-0}, thereby concluding the proof of the proposition.

\end{proof}

\section{Optimality remarks} 
\label{s:sharp}
As noted earlier, a key obstruction to global dispersive estimates in the focusing case
is given by the potential existence of small solitons. Here we will restrict our attention to the simplest model, namely the focusing cubic NLS, and  test the optimality of our estimates 
on solitons for this model. It is not so difficult to show that small solitons exist for our model whenever the focusing assumption (H4) is satisfied in some frequency region.

All cubic NLS solitons are equivalent modulo scaling and Galilean transformations. 
Our bounds are Galilean invariant, so we set the soliton velocity to zero and we
focus on scaling. Then the unit scale soliton has the form
\[
u(x,t) = e^{it}  Q(x), \qquad Q(x) = \sech (x).
\]
Rescaled to the frequency scale $\lambda$, this yields the solitons
\[
u_\lambda(x,t) = e^{it \lambda^2} Q_\lambda(x), 
\qquad Q_\lambda(x) = \lambda Q(\lambda x),
\]
which has initial data size 
\[
\|u_0\|_{L_x^2}^2 = \|Q_\lambda\|_{L^2_x}^2 
=\int_{\mathbb{R}}\lambda^2 Q^2( \lambda x)\, dx
\approx \lambda,
\]
so for Theorem~\ref{t:focusing} we will choose $\lambda = \epsilon^2$.

On the other hand,
\[
\Vert u_{\lambda}\Vert_{L^6_{t,x}(0, T, \mathbb{R})}^6 =T\int_{\mathbb{R}}\lambda^6  Q^6( \lambda x)\, dx\approx  T\lambda^5 .
\]
Then it is easily seen that we have approximate equality in \eqref{main-Str} and \eqref{main-Str-model}.

Similarly, we compute the bilinear Strichartz norm. Due to the bound from below 
on $|I|$, we have $c \gtrsim \lambda^{-\frac12}$ whereas $u$ is concentrated at frequency $\lesssim \lambda$. 
Then 
\[
\Vert \partial_x |u_{\lambda}|^2\Vert_{L^2_t(0, T; \dot H^{\frac12}_x + c L^2_x)}^2 
\approx c^{-2} \|  \partial_x |u_{\lambda}|^2\Vert_{L^2_{t,x}(0, T;  L^2)}^2 
=  \lambda (T \lambda^2)^{-1} T\int_{\mathbb{R}}\lambda^6 \partial_x Q_x^2( \lambda x)\, dx\approx  \lambda^2 .
\]
This corresponds to having equality in \eqref{main-bi} and \eqref{main-bi-model}.


\begin{thebibliography}{10}

\bibitem{AD}
Thomas Alazard and Jean-Marc Delort.
\newblock Global solutions and asymptotic behavior for two dimensional gravity
  water waves.
\newblock {\em Ann. Sci. \'{E}c. Norm. Sup\'{e}r. (4)}, 48(5):1149--1238, 2015.

\bibitem{IST-focusing}
Michael Borghese, Robert Jenkins, and Kenneth D. T.-R. McLaughlin.
\newblock Long time asymptotic behavior of the focusing nonlinear
  {S}chr\"{o}dinger equation.
\newblock {\em Ann. Inst. H. Poincar\'{e} C Anal. Non Lin\'{e}aire},
  35(4):887--920, 2018.

\bibitem{Bourgain-nf}
Jean Bourgain.
\newblock A remark on normal forms and the ``{$I$}-method'' for periodic {NLS}.
\newblock {\em J. Anal. Math.}, 94:125--157, 2004.

\bibitem{MR2527809}
J.~Colliander, M.~Grillakis, and N.~Tzirakis.
\newblock Tensor products and correlation estimates with applications to
  nonlinear {S}chr\"{o}dinger equations.
\newblock {\em Comm. Pure Appl. Math.}, 62(7):920--968, 2009.

\bibitem{I-method}
J.~Colliander, M.~Keel, G.~Staffilani, H.~Takaoka, and T.~Tao.
\newblock Almost conservation laws and global rough solutions to a nonlinear
  {S}chr\"{o}dinger equation.
\newblock {\em Math. Res. Lett.}, 9(5-6):659--682, 2002.

\bibitem{MR2053757}
J.~Colliander, M.~Keel, G.~Staffilani, H.~Takaoka, and T.~Tao.
\newblock Global existence and scattering for rough solutions of a nonlinear
  {S}chr\"{o}dinger equation on {$\mathbb R^3$}.
\newblock {\em Comm. Pure Appl. Math.}, 57(8):987--1014, 2004.

\bibitem{MR2415387}
J.~Colliander, M.~Keel, G.~Staffilani, H.~Takaoka, and T.~Tao.
\newblock Global well-posedness and scattering for the energy-critical
  nonlinear {S}chr\"{o}dinger equation in {$\mathbb R^3$}.
\newblock {\em Ann. of Math. (2)}, 167(3):767--865, 2008.

\bibitem{I-method2}
J.~Colliander, M.~Keel, G.~Staffilani, H.~Takaoka, and T.~Tao.
\newblock Resonant decompositions and the {$I$}-method for the cubic nonlinear
  {S}chr\"{o}dinger equation on {$\mathbb R^2$}.
\newblock {\em Discrete Contin. Dyn. Syst.}, 21(3):665--686, 2008.

\bibitem{IST}
Percy Deift and Xin Zhou.
\newblock Long-time asymptotics for solutions of the {NLS} equation with
  initial data in a weighted {S}obolev space.
\newblock {\em Comm. Pure Appl. Math.}, 56(8):1029--1077, 2003.
\newblock Dedicated to the memory of J\"{u}rgen K. Moser.

\bibitem{D}
Jean-Marc Delort.
\newblock Semiclassical microlocal normal forms and global solutions of
  modified one-dimensional {KG} equations.
\newblock {\em Ann. Inst. Fourier (Grenoble)}, 66(4):1451--1528, 2016.

\bibitem{MR3483476}
Benjamin Dodson.
\newblock Global well-posedness and scattering for the defocusing, {$L^2$}
  critical, nonlinear {S}chr\"{o}dinger equation when {$d=1$}.
\newblock {\em Amer. J. Math.}, 138(2):531--569, 2016.

\bibitem{MR3625190}
Benjamin Dodson.
\newblock Global well-posedness and scattering for the defocusing,
  mass-critical generalized {K}d{V} equation.
\newblock {\em Ann. PDE}, 3(1):Paper No. 5, 35, 2017.

\bibitem{HN}
Nakao Hayashi and Pavel~I. Naumkin.
\newblock Asymptotics for large time of solutions to the nonlinear
  {S}chr\"{o}dinger and {H}artree equations.
\newblock {\em Amer. J. Math.}, 120(2):369--389, 1998.

\bibitem{HN1}
Nakao Hayashi and Pavel~I. Naumkin.
\newblock Large time asymptotics for the fractional nonlinear {S}chr\"{o}dinger
  equation.
\newblock {\em Adv. Differential Equations}, 25(1-2):31--80, 2020.

\bibitem{hi}
John~K. Hunter and Mihaela Ifrim.
\newblock Enhanced life span of smooth solutions of a {B}urgers-{H}ilbert
  equation.
\newblock {\em SIAM J. Math. Anal.}, 44(3):2039--2052, 2012.

\bibitem{BH}
John~K. Hunter, Mihaela Ifrim, Daniel Tataru, and Tak~Kwong Wong.
\newblock Long time solutions for a {B}urgers-{H}ilbert equation via a modified
  energy method.
\newblock {\em Proc. Amer. Math. Soc.}, 143(8):3407--3412, 2015.

\bibitem{IT-NLS}
Mihaela Ifrim and Daniel Tataru.
\newblock Global bounds for the cubic nonlinear {S}chr\"{o}dinger equation
  ({NLS}) in one space dimension.
\newblock {\em Nonlinearity}, 28(8):2661--2675, 2015.

\bibitem{IT-g}
Mihaela Ifrim and Daniel Tataru.
\newblock Two dimensional water waves in holomorphic coordinates {II}: {G}lobal
  solutions.
\newblock {\em Bull. Soc. Math. France}, 144(2):369--394, 2016.

\bibitem{IT-c}
Mihaela Ifrim and Daniel Tataru.
\newblock The lifespan of small data solutions in two dimensional capillary
  water waves.
\newblock {\em Arch. Ration. Mech. Anal.}, 225(3):1279--1346, 2017.

\bibitem{IT-BO}
Mihaela Ifrim and Daniel Tataru.
\newblock Well-posedness and dispersive decay of small data solutions for the
  {B}enjamin-{O}no equation.
\newblock {\em Ann. Sci. \'{E}c. Norm. Sup\'{e}r. (4)}, 52(2):297--335, 2019.

\bibitem{IT-noloc}
Mihaela {Ifrim} and Daniel {Tataru}.
\newblock {Global solutions for 1D cubic defocusing dispersive equations: Part
  I}.
\newblock {\em arXiv e-prints}, page arXiv:2205.12212, May 2022.

\bibitem{IT-wp}
Mihaela {Ifrim} and Daniel {Tataru}.
\newblock {Testing by wave packets and modified scattering in nonlinear
  dispersive pde's}.
\newblock {\em arXiv e-prints}, page arXiv:2204.13285, April 2022.

\bibitem{KP}
Jun Kato and Fabio Pusateri.
\newblock A new proof of long-range scattering for critical nonlinear
  {S}chr\"{o}dinger equations.
\newblock {\em Differential Integral Equations}, 24(9-10):923--940, 2011.

\bibitem{LLS}
Hans Lindblad, Jonas L\"{u}hrmann, and Avy Soffer.
\newblock Asymptotics for 1{D} {K}lein-{G}ordon equations with variable
  coefficient quadratic nonlinearities.
\newblock {\em Arch. Ration. Mech. Anal.}, 241(3):1459--1527, 2021.

\bibitem{LS}
Hans Lindblad and Avy Soffer.
\newblock Scattering and small data completeness for the critical nonlinear
  {S}chr\"{o}dinger equation.
\newblock {\em Nonlinearity}, 19(2):345--353, 2006.

\bibitem{PV}
Fabrice Planchon and Luis Vega.
\newblock Bilinear virial identities and applications.
\newblock {\em Ann. Sci. \'{E}c. Norm. Sup\'{e}r. (4)}, 42(2):261--290, 2009.

\bibitem{MR2288737}
E.~Ryckman and M.~Visan.
\newblock Global well-posedness and scattering for the defocusing
  energy-critical nonlinear {S}chr\"{o}dinger equation in {$\mathbb R^{1+4}$}.
\newblock {\em Amer. J. Math.}, 129(1):1--60, 2007.

\bibitem{Shatah}
Jalal Shatah.
\newblock Normal forms and quadratic nonlinear {K}lein-{G}ordon equations.
\newblock {\em Comm. Pure Appl. Math.}, 38(5):685--696, 1985.

\bibitem{Tao-WM}
Terence Tao.
\newblock Global regularity of wave maps. {II}. {S}mall energy in two
  dimensions.
\newblock {\em Comm. Math. Phys.}, 224(2):443--544, 2001.

\bibitem{Tao-BO}
Terence Tao.
\newblock Global well-posedness of the {B}enjamin-{O}no equation in {$H^1({\bf
  R})$}.
\newblock {\em J. Hyperbolic Differ. Equ.}, 1(1):27--49, 2004.

\end{thebibliography}

\bibliographystyle{plain}

\end{document}